\documentclass[final,3p,times]{elsarticle}

\usepackage[normalem]{ulem}
\usepackage{graphicx}
\usepackage{amsmath}
\usepackage{mathtools}
\usepackage{amsfonts}
\usepackage{amssymb}
\usepackage{amsthm}

\newcommand{\R}{ \mathbb{R}}
\newcommand{\Z}{ \mathbb{Z}}

\newcommand{\Lop}{{\rm L}}

\newcommand{\dint}{{\rm d}}

\newcommand{\Fourier}{ \mathcal{F}}

\newcommand{\x}{r}

\newcommand{\eop}{\hfill$\Box$}

\newenvironment{pf}%
{\par\noindent\textbf{Proof:}~}%

\newtheorem{definition}{Definition}
\newtheorem{proposition}{Proposition}

\newtheorem{theorem}{Theorem}
\newtheorem{remark}{Remark}

\def\M#1{{\bf{#1}}}  

\def\ee{\mathrm{e}}     

\def\dint{\;\mathrm{d}}

\begin{document}

\begin{frontmatter}

\title{Ellipse-Preserving Hermite interpolation and Subdivision}

\author[label1]{Costanza Conti}
\ead{costanza.conti@unifi.it}

\author[label2]{Lucia Romani}
\ead{lucia.romani@unimib.it}

\author[label3]{Michael Unser}
\ead{Michael.Unser@epfl.ch}

\address[label1]{Dipartimento di Energetica ``Sergio Stecco'',
Universit\`{a} di Firenze,\\
Viale Morgagni 40/44, 50134 Firenze, Italy\\
Tel.: +39-0554796713,
Fax: +39-0554224137}

\address[label2]{Department of Mathematics and Applications, University of Milano-Bicocca,\\
Via R. Cozzi 53, 20125 Milano, Italy\\
Tel.: +39-0264485735,
Fax: +39-0264485705}

\address[label3]{Biomedical Imaging Group,
Ecole Polytechnique F\'{e}d\'{e}rale de Lausanne,\\
LIB, BM 4.136 (Batiment BM), Station 17, CH-1015 Lausanne, Switzerland\\
Tel.: +41(21)693.51.75,
Fax : +41(21)693.37.01}

\begin{abstract}
We introduce a family of piecewise-exponential functions that have the Hermite interpolation property. Our design is motivated by the search for an effective scheme for the joint interpolation of points and associated tangents on a curve with the ability to perfectly reproduce ellipses. We prove that the proposed Hermite functions form a Riesz basis and that they reproduce prescribed exponential polynomials.
We present a method based on Green's functions to unravel their multi-resolution and approximation-theoretic properties.
Finally, we derive the corresponding vector and scalar subdivision schemes, which lend themselves to a fast implementation.
The proposed vector scheme is interpolatory and level-dependent, but its asymptotic behaviour is the same as the classical cubic Hermite spline algorithm. The same convergence properties---i.e., fourth order of approximation---are hence ensured.
\end{abstract}

%
%

\end{frontmatter}

\section{Introduction}

Cubic Hermite splines are piecewise-cubic polynomial functions that are parametrized in terms of the value of the function and its derivative at the end point of each polynomial segment. By construction, the resulting spline is continuous with continuous first-order derivative. Cubic Hermite splines are used extensively in computer graphics and geometric modelling to represent curves as well as motion trajectories that pass through specified anchor points with prescribed tangents \cite{HL93}. This is typically achieved by fitting a separate Hermite spline interpolant for each coordinate variable.

Cubic Hermite splines have a number of attractive computational features. The basis functions are interpolating with a fourth-order approximation and their support is minimal.
They satisfy multiresolution properties, which is the key to the specification of  subdivision schemes \cite{Merrien} and the construction of multi-wavelet bases \cite{DaHanJiaKu00,Plonka95}. They are also closely linked to the B\'{e}zier curves, which provide an equivalent mode of representation.
Their only limitation is that they require many control points to accurately reproduce elementary shapes such as circles and ellipses. This is why we investigate in this paper a variation of the classical Hermite scheme that is specifically geared towards the reproduction of elliptical shapes. These exponential Hermite splines are  ideally suited for outlining roundish objects in images with few control points (see \cite{Noi2014} for an application of this model to the segmentation of biomedical images).
Our main point in this work will be to show that we are able to achieve perfect ellipse reproduction while retaining all the attractive properties of the cubic Hermite splines modulo some proper adjustment of the underlying computational machinery. The extended Hermite functions that we shall specify are piecewise exponential polynomials with pieces in ${\cal E}_4:=\langle 1,x,\ee^{i\omega_0 x},\ee^{-i\omega_0 x}\rangle$, $\omega_0\in [0,\pi]$, joining $C^1$-continuously at the integer knots.
Hence they belong to the class of trigonometric splines.
This points towards a connection with other exponential spline basis functions investigated in the literature
(see, e.g., \cite{Dahmen:Micchelli:1987,MainarPena2007,ManniMazure10,Mazure1,Mazure05,McCartin:1991} and references quoted therein), although we are not aware of any prior work that specifically addresses the problem of ellipse reproduction nor covers the theoretical results that we are reporting here.

\smallskip \noindent The paper is organised as follows. In Section \ref{sec:2}, we motivate our design while spelling out the conditions that the basis functions must satisfy. We then derive the two  Hermite functions $(\phi_{1,\omega_0}, \phi_{2,\omega_0})$ that fulfil our requirements in Section \ref{sec:3}; these are the generators for the space $S_{{\cal E}_4}^1(\Z)$, which is made up of functions that are piecewise exponential polynomials with (double) knots on the integers. In Section \ref{sec:4}, we make the connection with exponential splines explicit by expressing the generators in terms of the Green's functions of the differential operators $\Lop_{1,\omega_0}:=\frac{\dint^4}{\dint x^4}+\omega_0^2\frac{\dint^2}{\dint x^2}$ and $\Lop_{2,\omega_0}:=\frac{\dint^3}{\dint x^3}+\omega_0^2\frac{\dint}{\dint x}$ (whose corresponding E-spline spaces are denoted as  $S_{{\cal E}_4}(\Z)$ and $S_{{\cal E}_3}(\Z)$). In Section \ref{sec:5},
we prove that the integer translates form a Riesz basis by analyzing the corresponding Gramian matrix.  Section \ref{sec:6} is devoted to the characterisation of the Hermite representation on $h\Z$, while Section \ref{sec:7} focuses on the investigation of its multi-resolution properties and the derivation of the corresponding subdivision scheme. In Section \ref{sec:8}, we show that our exponential Hermite splines, in direct analogy with their polynomial counterpart, admit a B\'{e}zier representation that involves an exponential generalization of the classical Bernstein polynomials. Finally, in Section \ref{sec:9}, we exploit the B\'{e}zier connection to derive the exponential version of the four point scalar subdivision scheme for the classical Hermite splines \cite{Juttler}.

\section{Motivation for the construction}\label{sec:2}
An active contour (a.k.a.\ snake) is a computational tool for detecting and outlining objects in digital images. Its central component is a closed parametric curve that evolves spatially towards the contour of a target by minimizing a suitable energy functional \cite{McInerney96}. The most commonly-used curve models rely on B-spline basis functions \cite{Brigger2000}.

\begin{figure}
\centering
\includegraphics[width=7cm]{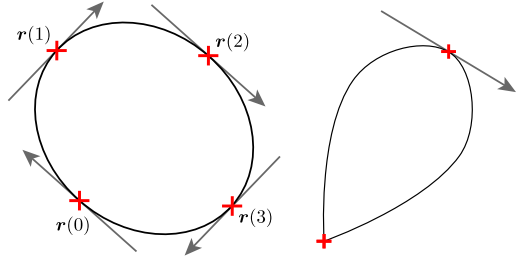}
\caption{\label{Fig:Snakes} Examples of parametric curves $\M r(t)$ represented in the Hermite basis.
The shape parameters  are the control points represented by crosses and the tangent handles (arrows) that control the derivative of each of the coordinate variable with respect to $t$. The first illustrates the ellipse-reproduction capability of our extended model, while the second demonstrates the production of a cusp by decreasing the magnitude of its tangent vector to zero.}
\end{figure}

Since roundish objects are common place in biological imaging (in particular, fluorescence microscopy), it is of interest
to develop a parametric framework that is specifically tailored to this type of shape while retaining the flexibility of splines and the ability to reparametrize by increasing the number of control points. A first solution to this problem was proposed by Delgado et al. who developed an ``active cells" framework that is based on cardinal exponential B-splines \cite{Del2012}.
The present research was motivated by the desire to refine this model by providing additional control over the tangents of the curve. This led us to the definition of a new parametric model that has the ability to perfectly reproducing ellipses while offering full tangential control as well as easy manipulation via the use of $M$ control points and B\'ezier handles. By introducing B\'ezier handles, one also gains in flexibility; for instance, one can induce a sharp break via a proper adjustment of the tangent vector (see Figure \ref{Fig:Snakes}).
The corresponding parametric representation is

\begin{align}
\label{eq:curve}
\M r(t)=\sum_{n\in \Z} \left(\M r(n)\phi_1(t-n) +  \M r'(n)\phi_2(t-n)\right)
\end{align}
where the closed curve
$\M r(t)=\big( x(t), y(t)\big)$ and its tangent $\M r'(t)=\big( \frac{\dint x(t)}{\dint t}, \frac{\dint y(t)}{\dint t}\big)$ are assumed to be $M$-periodic. Practically, this means that the underlying curve is uniquely specified by its shape parameters $\{\M r(n),\M r'(n)\}_{n=0}^{M-1}$ which can be translated graphically into a set of control points with tangential handles (see Figure \ref{Fig:Snakes}).

\begin{figure}
\centering
\includegraphics[width=6cm]{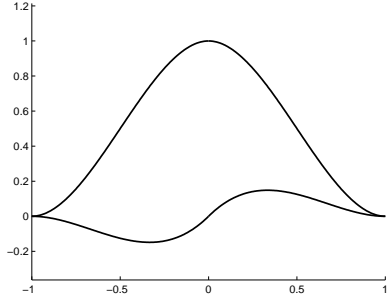}
\caption{\label{Fig:Hermite}The generators $\phi_{1,\omega_0}$ and $\phi_{2,\omega_0}$ of ${\cal E}_4$ Hermite splines with
$\omega_0=3/4\pi$. The two functions and their derivatives are vanishing at the integers with the exception of $\phi_{1,\omega_0}(0)=1$ and $\phi_{2,\omega_0}'(0)=1$ (interpolation conditions). Their support size is two.}
\end{figure}

The fundamental property of this kind of Hermite representation is that the generating functions $\phi_1$, $\phi_2$ and their derivatives $\phi'_1, \phi'_2$ satisfy the joint interpolation conditions
\begin{align*}
\phi_1(n)=\delta_{n,0}, \quad\phi'_2(n)=\delta_{n,0}, \quad
\phi'_1(n)=0, \quad
\phi_2(n)=0\end{align*}
for all $n \in \Z$ (see Figure \ref{Fig:Hermite}).

The shape space associated with \eqref{eq:curve} is the collection of all possible curves that can be generated by varying the control parameters $\{\M r(n),\M r'(n)\}_{n=0}^{M-1}$. Our three design requirements on the specification of this shape space are as follows:
\begin{enumerate}
\item the representation should be unambiguous and stable with respect to the variation of the shape parameters;
\item the shape space should be closed with respect to affine transformations;
\item the shape space should include all ellipses.
\end{enumerate}

The first point is taken care of by making sure that the basis functions form a Riesz basis  (see Section \ref{Sec:Riesz}).
We conclude this introductory section by providing the conditions on the basis functions that are imposed by the two latter requirements.

\subsection{Affine invariance}
Consider the affine transformation $\M s(t)=\M A\M r(t) + \M b$ of the curve $\M r(t)$ in the 2-D plane. We would like this new curve to be representable in the Hermite basis as
$$
\M A\M r(t) + \M b= \sum_{n\in \Z} \big(\underbrace{(\M A\M r(n)+\M b)}_{\M s(n)}\;\phi_1(t-n) + \underbrace{\M A \M r'(n)}_{\M s'(n)}\;\phi_2(t-n)\big),
$$
which is possible if and only if $\phi_1$ satisfies the partition of unity property
$$
\sum_{n\in \Z} \phi_1(t-n)= 1.
$$
\subsection{Reproduction of ellipses}
The specificity of our design is the ability to reproduce ellipses,
as illustrated in Figure \ref{Fig:Snakes}. Since the representation is affine invariant, it is sufficient to be able to encode the unit circle, which translates into the two complementary conditions
\begin{align*}
\cos( \omega_0 t)&=\sum_{n\in \Z} \big(\cos(\omega_0 n) \phi_1(t-n) -
\omega_0 \sin(\omega_0 n) \phi_2(t-n)\big)
\\
\sin(\omega_0 t)&=\sum_{n\in \Z} \big(\sin( \omega_0 n) \phi_1(t-n) +
 \omega_0 \cos( \omega_0 n) \phi_2(t-n)\big)
\end{align*}
with $\omega_0=\frac{2 \pi}{N} \in [0,\pi]$.

\section{Cardinal Hermite exponential splines}\label{sec:3}
In analogy with the classical cubic solution, we shall determine our extended Hermite functions $\phi_{1,\omega_0}(x)$ and $\phi_{2,\omega_0}(x)$ by first focusing on the unit interval $x \in[0,1]$ and imposing the four required boundary conditions in each case;
i.e., $$\phi_{1,\omega_0}(0)=1, \ \phi'_{1,\omega_0}(0)=0, \ \phi_{1,\omega_0}'(1)=0, \ \phi_{1,\omega_0}(1)=0$$
and $$\phi_{2,\omega_0}(0)=0, \ \phi'_{2,\omega_0}(0)=1, \ \phi_{2,\omega_0}'(1)=0, \ \phi_{2,\omega_0}'(1)=0.$$
The existence of such functions is guaranteed if we consider a common four-dimensional solution space of Tchebychev polynomials.
Because of our reproduction requirements, we already know that the solution space should contain the functions $\{1,\cos( \omega_0 x), \sin( \omega_0 x)\}$.
The last functional degree of freedom is taken care of by imposing that the two generators, which are supported in $[-1,1]$, should be real-valued, symmetric or anti-symmetric and restricted to the class of exponential polynomials in order to yield {\it bona fide} splines.
This fixes the solution space to  ${\cal E}_4=\langle 1,\ee^{i \omega_0x}, \ee^{-i \omega_0x},x \rangle$ and makes the construction of our trigonometric splines unique.

The functions $\phi_{1,\omega_0}$ and $\phi_{2,\omega_0}$ that fulfill these constraints then constitute the generators for the space of  cardinal Hermite exponential splines which is denoted by $S_{{\cal E}_4}^1(\Z)$.
They are given
by
\begin{equation}\label{eq:phi12_glob}
\phi_{1,\omega_0}(x)= \left\{
\begin{array}{ll}
g_{1,\omega_0 }(x), & \mbox{for } x \geq 0\\
g_{1,\omega_0 }(-x), & \mbox{for } x<0
\end{array}
\right., \qquad
\phi_{2,\omega_0}(x)= \left\{
\begin{array}{ll}
g_{2,\omega_0 }(x), & \mbox{for } x \geq 0\\
-g_{2,\omega_0 }(-x), & \mbox{for } x<0
\end{array}
\right.
\end{equation}
where
\begin{equation}\label{def:g}
g_{j,\omega_0 }(x):=\left(a_j(\omega_0) + b_j(\omega_0) x + c_j(\omega_0) \ee^{i \omega_0 x} + d_j(\omega_0) \ee^{-i\omega_0 x}\right)\chi_{[0,1]}, \quad j=1,2
\end{equation}

\smallskip \noindent
with the combination coefficients given by
$$
\begin{array}{l}
a_1(\omega_0) := \frac{i\omega_0+1 + \ee^{i\omega_0} (i\omega_0-1)}{q(\omega_0)}, \quad b_1(\omega_0) := -\frac{i\omega_0 (\ee^{i\omega_0} + 1)}{q(\omega_0)}, \quad
c_1(\omega_0) := \frac{1}{q(\omega_0)}, \quad d_1(\omega_0) := -\frac{\ee^{i\omega_0}}{q(\omega_0)},
\end{array}
$$
$$
\begin{array}{l}
a_2(\omega_0):= \frac{p(\omega_0)}{i\omega_0 (\ee^{i\omega_0}-1) q(\omega_0) },\quad b_2(\omega_0):= -\frac{\ee^{i\omega_0}-1}{q(\omega_0)},\quad
c_2(\omega_0):= \frac{\ee^{i\omega_0}-i\omega_0-1}{i\omega_0(\ee^{i\omega_0}-1)q(\omega_0)},\quad
d_2(\omega_0):= -\frac{\ee^{i\omega_0} (\ee^{i\omega_0}(i\omega_0-1) + 1)}{i\omega_0(\ee^{i\omega_0} - 1) q(\omega_0)},\\
\end{array}
$$

\medskip \noindent
and
\begin{equation}\label{def:p&q}
p(\omega_0):=\ee^{2i\omega_0}(i\omega_0-1)+i\omega_0+1,\qquad q(\omega_0):=\ee^{i\omega_0}(i\omega_0-2)+i\omega_0+2.
\end{equation}

\medskip  \noindent
Since $a_j(\omega_0),\ b_j(\omega_0),\ j=1,2$ are both real as well as $c_j(\omega_0) \ee^{i \omega_0 x} + d_j(\omega_0) \ee^{-i\omega_0 x},\ j=1,2$, both functions in \eqref{def:g} are real-valued. Indeed substitution of the above coefficients in (\ref{def:g}) provides
\begin{equation}\label{def:g_real-version}
\begin{array}{ll}
g_{1,\omega_0 }(x)&=\left(1-\frac{\sin(\omega_0/2)}{s(\omega_0)} + \frac{\omega_0\cos(\omega_0/2)}{s(\omega_0)} x + \frac{\sin(\omega_0/2-\omega_0x)}{s(\omega_0)}\right)\chi_{[0,1]},\\
\\
g_{2,\omega_0 }(x)&=\left(\frac{\sin(\omega_0)-\omega_0\cos(\omega_0)}{\omega_0 \, u(\omega_0)} +
 \frac{\sin(\omega_0/2)}{s(\omega_0)}x- \frac{\omega_0^2 \cos(\omega_0/2) \cos(\omega_0(1-x)) + \sin(\omega_0/2) \big( \sin(\omega_0 x) u(\omega_0) -\cos(\omega_0 x) v(\omega_0) \big) }{2\omega_0 \sin(\omega_0/2) s(\omega_0) \, t(\omega_0)} \right)\chi_{[0,1]},
\end{array}
\end{equation}
where
\begin{equation}\label{eq:sw0}
s(\omega_0):=2\sin(\omega_0/2)-\omega_0\cos(\omega_0/2), \qquad  t(\omega_0):=2\sin(\omega_0/2)+\omega_0\cos(\omega_0/2),
\end{equation}
and
\begin{equation}\label{eq:u&v}
u(\omega_0):=\omega_0 \sin(\omega_0)-2 (1- \cos(\omega_0)), \qquad v(\omega_0):=2 \sin(\omega_0)+\omega_0 (1- \cos(\omega_0)).
\end{equation}

\smallskip  \noindent
Note that $\phi_{1,\omega_0}$ and $\phi_{2,\omega_0}$ are exponential polynomials in ${\cal E}_4$ in each interval $[n,n+1)$ for $n=-1,0$ (and by extension for any $n \in \Z$) and that they are differentiable (with continuous derivatives) at the knots $x=n$.

\smallskip  \noindent
It is clear that any linear combinations of the integer shifts of these functions is a piecewise exponential polynomial made of pieces in ${\cal E}_4$ joining $C^1$-continuously at the integers. Such functions can also be interpreted as exponential splines with double knots on the integers, the effect of a double knot being to reduce
the ordinary degree of continuity of the classical cardinal exponential splines by one \cite{Kunkle99}.
It follows that the space $S_{{\cal E}_4}^1(\Z)$ can be written as
\begin{align}
\label{eq:S41}
S_{{\cal E}_4}^1(\Z)=\left\{s(x)=\sum_{n \in \Z} \M a^T[n]\, \boldsymbol\phi_{\omega_0}(x-n): \M a \in \ell_2^{2\times 1}(\Z) \right\}, \qquad \hbox{with} \quad \boldsymbol \phi_{\omega_0}:=(\phi_{1,\omega_0},\phi_{2,\omega_0})^T.
\end{align}
Due to the Hermite interpolation condition, the expansion coefficients in \eqref{eq:S41}
coincide with the samples of the function and its first derivative on the integer grid; that is, $\M a[n]=\big(s(n),s'(n)\big)^T$.

\section{Connection with standard exponential splines and reproduction properties}\label{sec:4}
Our way of establishing the link with standard exponential splines is to compute the Fourier transforms of the Hermite exponential spline generators
with the convention that $\hat f(\omega)=\int_\R f(x) \ee^{-i \omega x} \dint x$.
This yields
\begin{equation}\label{eq:phihat}
\widehat{\boldsymbol \phi}_{\omega_0}(\omega)=
\left[
\begin{array}{c}
\hat \phi_{1,\omega_0}(\omega) \\
\hat \phi_{2,\omega_0}(\omega)
\end{array}
\right]=
\left[
\begin{array}{c}
\frac{c_0^1(\omega_0)+c_1^1(\omega_0)\omega}{\omega^2(\omega^2-\omega_0^2)} \smallskip\\
\frac{c_0^2(\omega_0)+c_1^2(\omega_0)\omega}{\omega^2(\omega^2-\omega_0^2)}
\end{array}
\right],
\end{equation}
where, for $p(\omega_0)$, $q(\omega_0)$ as in \eqref{def:p&q}, we have
$$
\begin{array}{ll}
c_0^1(\omega_0):= \frac{-i\omega_0^3(\ee^{i\omega_0}+1)}{q(\omega_0)} (2-(\ee^{-i\omega} + \ee^{i\omega})), &  c_1^1(\omega_0):=  \frac{-i \omega_0^2 (\ee^{i\omega_0}-1)}{q(\omega_0)} (\ee^{i\omega}-\ee^{-i\omega}),\\
c_0^2(\omega_0):= \frac{-\omega_0^2(\ee^{i\omega_0}-1)}{q(\omega_0)} (\ee^{i\omega}-\ee^{-i\omega}),
&  c_1^2(\omega_0):= \frac{- \omega_0}{q(\omega_0)(1-\ee^{i\omega_0})} \Big(2p(\omega_0)-(1+2 i \omega_0 \ee^{i\omega_0}-\ee^{2i\omega_0})(\ee^{-i\omega}+\ee^{i\omega})\Big).
\end{array}
$$
Next, we rewrite (\ref{eq:phihat}) in matrix-vector form as
\begin{align}
\label{Eq:FouriertoGreen}
\widehat{\boldsymbol \phi}_{\omega_0}(\omega) =\widehat{\M R}(\ee^{i \omega }) \ \widehat{\boldsymbol \rho}_{\omega_0}(\omega),
\end{align}
with
\begin{equation}\label{eq:rhohat}
\widehat{\boldsymbol \rho}_{\omega_0}(\omega):=\left[
\begin{array}{c}
\hat \rho_{1,\omega_0}(\omega) \\
\hat \rho_{2,\omega_0}(\omega)
\end{array}
\right]=\left[
\begin{array}{c}
\frac{1}{\omega^2(\omega^2-\omega_0^2)} \\
\frac{i}{\omega (\omega^2-\omega_0^2)}
\end{array}
\right],
\end{equation}
and
\begin{eqnarray}\label{eq:Rmat}
\widehat{\M R}(\ee^{i \omega })&:=&\frac{i\omega_0}{q(\omega_0)}
\left[
\begin{array}{cc}
 -\omega_0^2 (\ee^{i\omega_0}+1) \Big(2-(\ee^{-i\omega}+\ee^{i\omega})\Big) & i\omega_0(\ee^{i\omega_0}-1)(\ee^{i\omega}-\ee^{-i\omega}) \smallskip\\
 i\omega_0(\ee^{i\omega_0}-1) (\ee^{i\omega}-\ee^{-i\omega}) & \frac{2p(\omega_0)-(1+2i\omega_0 \ee^{i\omega_0}-\ee^{2i \omega_0})(\ee^{-i\omega}+\ee^{i\omega})}{1-\ee^{i\omega_0}}
\end{array}
\right] \medskip\\
&=& \frac{1}{s(\omega_0)} \left[
\begin{array}{cc}
\omega_0^3 \cos(\omega_0/2) \Big(2-(\ee^{-i\omega}+\ee^{i\omega})\Big) &  \omega_0^2 \sin(\omega_0/2) (\ee^{i\omega}-\ee^{-i\omega}) \smallskip\\
\omega_0^2 \sin(\omega_0/2) (\ee^{i\omega}-\ee^{-i\omega}) &  \frac{2 \omega_0 \big(\omega_0 \cos(\omega_0)-\sin(\omega_0)\big)-\omega_0 \big(\omega_0-\sin(\omega_0)\big) (\ee^{-i\omega}+\ee^{i\omega})}{2 \sin(\omega_0/2)}
\end{array}
\right], \nonumber
\end{eqnarray}
where expressions of $p(\omega_0)$, $q(\omega_0)$ and $s(\omega_0)$ are given in (\ref{def:p&q}) and (\ref{eq:sw0}), respectively.
The aim here is to reveal the linear relation between the
generators ${\boldsymbol \phi}_{\omega_0}=(\phi_{1,\omega_0},\phi_{2,\omega_0})^T$ and  ${\boldsymbol \rho}_{\omega_0}=(\rho_{1,\omega_0},\rho_{2,\omega_0})^T$. The latter 
are the Green's
functions of the differential operators $\Lop_{1,\omega_0}:=\frac{\dint^4}{\dint x^4}+\omega_0^2\frac{\dint^2}{\dint x^2}$, $\Lop_{2,\omega_0}:=\frac{\dint^3}{\dint x^3}+\omega_0^2\frac{\dint}{\dint x}$ defining $S_{{\cal E}_4}(\Z)$ and $S_{{\cal E}_3}(\Z)$, respectively. The explicit expression of these Green's functions is given by
\begin{align}
\rho_{1,\omega_0}(x)&=\Fourier^{-1}\left\{ \frac{1}{\omega^2(\omega^2-\omega_0^2)}\right\}=
\frac{{\omega_0} x-\sin({\omega_0} x)}{2{\omega_0}^3} \, \text{sgn}(x),
\\
\rho_{2,\omega_0}(x)&=\frac{\dint}{\dint x} \rho_{1,\omega_0}(x)=
\frac{1-\cos({\omega_0} x)}{2{\omega_0}^2} \, \text{sgn}(x).
\end{align}

\noindent By inverting the $2 \times 2$ Fourier matrix $\widehat{\M R}(\ee^{i \omega })$ in (\ref{eq:Rmat}), we find that
\begin{align}
\widehat{\boldsymbol \rho}_{\omega_0}(\omega)=\Big(\widehat{\M R}(\ee^{i \omega })\Big)^{-1} \ \widehat{\boldsymbol \phi}_{\omega_0}(\omega).
\end{align}
Since $\Big(\widehat{\M R}(\ee^{i \omega })\Big)^{-1}=:\widehat{\M P}(\ee^{i \omega })$ has entries that are ratios of trigonometric polynomials, its discrete-time inverse Fourier transform is well-defined and guaranteed to yield a unique sequence of matrices
$$\M P[n]=\frac{1}{2 \pi}\int_{-\pi}^{+\pi} \widehat{\M P}(\ee^{i \omega }) \ \ee^{i \omega n} \dint \omega$$
of slow growth. Hence, we conclude that
\begin{align}
\boldsymbol \rho_{\omega_0}(x)
=\sum_{n \in \Z} \M P[n] \,\boldsymbol \phi_{\omega_0}(x-n),
\label{eq:Greenrepro1}
\end{align}
which proves that the Green's functions $\rho_{1,\omega_0}$ and $\rho_{2,\omega_0}$ (as well as their integer shifts) can be perfectly reproduced by $\{\boldsymbol \phi_{\omega_0}(\cdot-n)\}_{n \in \Z}$.
The specific form of \eqref{eq:Greenrepro1} then follows from the interpolation property of the generators; that is, from the relation
\begin{align}
\label{Eq:Hermiteprop}
s(x)= \sum_{n \in \Z} \big(s(n)\phi_{1,\omega_0}(x-n)+s'(n)\phi_{2,\omega_0}(x-n)\big),
\end{align}
which is valid for any function in $S^1_{{\cal E}_4}(\Z)$.
In particular, we have that
\begin{equation}\label{def:ro_1versus-phi}
\rho_{1,\omega_0}(x)=\sum_{n \in \Z}  \left( \frac{{\omega_0} n-\sin({\omega_0} n)}{2{\omega_0}^3} \, \text{sgn}(n)  \phi_{1,\omega_0}(x-n) + \frac{1-\cos(\omega_0n)}{2{\omega_0}^2} \, \text{sgn}(n) \phi_{2,\omega_0}(x-n)\right)\,,
\end{equation}
and
\begin{equation}\label{def:ro_2versus-phi}
\rho_{2,\omega_0}(x)=\sum_{n \in \Z}\left( \frac{1-\cos(\omega_0n)}{2{\omega_0}^2} \, \text{sgn}(n) \phi_{1,\omega_0}(x-n) + \frac{\sin(\omega_0 n)}{2{\omega_0}} \, \text{sgn}(n) \phi_{2,\omega_0}(x-n)\right).
\end{equation}

\smallskip
In order to establish a link between order-four Hermite exponential splines and order-four exponential B-splines (see \cite{Unser2005} for the definition and detailed investigation of exponential B-splines), we consider a discretization on ${\mathbb Z}$ of the differential operators $\Lop_{1,\omega_0}$ and $\Lop_{2,\omega_0}$ based on the following recursive definition of what we call the
discrete annihilation operator. The basic principle here is to specify the shortest possible sequence of weights that annihilates the components (typically sinusoids) that are in the null space of those operators.
\begin{definition}\label{Delta}
For $\omega_j \in [0,\pi],\ j=0,\dots,m$, the discrete annihilation operator for the frequencies $(\omega_0,\cdots,\omega_m)$ is recursively defined as
\begin{equation}\label{def:delta}
\Delta_{\omega_0} f(x):=f(x)-\ee^{i\omega_0} f(x-1),\qquad \Delta_{(\omega_0,\cdots,\omega_m)} f(x):=\Delta_{\omega_0} \left(\Delta_{(\omega_1,\cdots,\omega_m)} f(x)\right)\,.
\end{equation}
\end{definition}

\smallskip \noindent
In light of the above definition, a discretization on ${\mathbb Z}$ of the differential operators $\Lop_{1,\omega_0}$ and $\Lop_{2,\omega_0}$ is given by $\Delta_{(0,0,\omega_0,-\omega_0)}$ and $\Delta_{(0,\omega_0,-\omega_0)}$, respectively. Note that $\Delta_{(0,0,\omega_0,-\omega_0)}$ is exact when applied to functions in ${\cal E}_4$ and that $\Delta_{(0,\omega_0,-\omega_0)}$ is exact when applied to functions in ${\cal E}_3:=\langle 1,\ee^{i\omega_0 x},\ee^{-i\omega_0 x}\rangle$, $\omega_0 \in [0,\pi]$.

\medskip \noindent In accordance with the classical theory of exponential splines, the order-four and order-three normalized exponential B-splines are defined as follows, with a normalization factor that ensures the partition of unity \cite{Christensen2012}.

\begin{definition}\label{def:Bexp43}
The normalized order-four exponential B-spline basis for $S_{{\cal E}_4}(\Z)$ is defined by
\begin{equation}\label{Eq:Bcubic_bis}
B_{4,\omega_0}(x)
=\left( \frac{\omega_0}{2 \sin(\omega_0/2)}\right)^2 \Delta_{(0,0,\omega_0,-\omega_0)} \, \rho_{1,\omega_0} (x)\,.
\end{equation}
Similarly, the normalized order-three exponential B-spline basis for $S_{{\cal E}_3}(\Z)$ is
\begin{equation}\label{Eq:Bquad_bis}
B_{3,\omega_0}(x)
=\left( \frac{\omega_0}{2 \sin(\omega_0/2)}\right)^2 \Delta_{(0,\omega_0,-\omega_0)} \, \rho_{2,\omega_0}(x).
\end{equation}
\end{definition}
\smallskip \noindent With the help of some algebra, we are also able to express $B_{4,\omega_0}$ and $B_{3,\omega_0}$ in terms of shifts of the generator $\boldsymbol\phi_{\omega_0}$. For instance, we find that
\begin{equation}
B_{4,\omega_0}(x)=\gamma_1^3 \phi_{1,\omega_0}(x-1) + \gamma_2^3 \phi_{1,\omega_0}(x - 2) + \gamma_3^3 \phi_{1,\omega_0}(x - 3) + \mu_1^3 \phi_{2,\omega_0}(x-1) +\mu_2^3 \phi_{2,\omega_0}(x-2) +
 \mu_3^3 \phi_{2,\omega_0}(x - 3),
 \label{Eq:Bcubic}
\end{equation}
where
$$
\gamma_1^3
=\frac{\omega_0 - \sin(\omega_0)}{4 \omega_0 \sin^2(\omega_0/2)}, \quad
\gamma_2^3=1-2\gamma_1^3,
\quad
\gamma_3^3=\gamma_1^3, \quad
\mu_1^3=\frac{1}{2}, \quad \mu_2^3=0, \quad \mu_3^3=-\frac{1}{2}.
$$
One can easily verify that $B_{4,\omega_0}$ is supported on $[0,4]$ and that it converges to a cubic B-spline as $\omega_0\to 0$.

Similarly, we make use of the Hermite interpolation property \eqref{Eq:Hermiteprop} to obtain the corresponding expression for the order-three exponential B-spline for the space $S_{{\cal E}_3}(\Z)$, which is, instead, supported on $[0,3]$:
\begin{equation}\label{def:phi2}
B_{3,\omega_0}(x)=\gamma_1^2 \phi_{1,\omega_0}(x-1) + \gamma_2^2 \phi_{1,\omega_0}(x - 2) +  \mu_1^2 \phi_{2,\omega_0}(x - 1) + \mu_2^2 \phi_{2,\omega_0}(x-2),
\end{equation}
where
$$
\gamma_1^2=\frac{1}{2}, \quad \gamma_2^2=\frac{1}{2}, \qquad
\mu_1^2
=\frac{\omega_0}{2} \cot(\omega_0/2), \quad \mu_2^2=-\mu_1^2.
$$
Since exponential B-splines reproduce functions in ${\cal E}_4$, the property automatically extends to the space $S_{{\cal E}_4}^1(\Z)$. Specifically, we have that
\begin{equation}\label{eq:pol_repro}
x^m= \sum_{n \in \Z} \big(n^m \phi_{1,\omega_0}(x-n)+ m\, n^{m-1} \phi_{2,\omega_0}(x-n)\big), \qquad m=0,1,
\end{equation}
and
\begin{equation}\label{eq:exp_repro}
\ee^{\pm i\omega_0 x}= \sum_{n \in \Z} \big(\ee^{\pm i\omega_0 n} \phi_{1,\omega_0}(x-n) \pm i\omega_0 \, \ee^{\pm i\omega_0 n} \phi_{2,\omega_0}(x-n)\big).
\end{equation}
\begin{remark}
From (\ref{eq:pol_repro}) and (\ref{eq:exp_repro}), we immediately observe that any Hermite interpolant of type (\ref{Eq:Hermiteprop}) is reproducing the whole space  ${\cal E}_4$ and, in particular, it is \emph{ellipse-reproducing}. Moreover, we observe that (\ref{Eq:Bcubic}) and (\ref{def:phi2}) can be interpreted as the construction of the shortest superfunction for the space $S_{{\cal E}_4}^1(\Z)$ \cite{CJ00}.
\end{remark}

\noindent The remarkable property with respect to the theory of exponential splines is that the space $S_{{\cal E}_4}^1(\Z)$,
(which is the sum of $S_{{\cal E}_4}(\Z)$ and $S_{{\cal E}_3}(\Z)$ as shown below), admits basis functions of size $2$ that are shorter than the exponential B-splines for any of the pure spline constituents. This can be explained via the so-called {\em localisation} process.
Based on \eqref{eq:phihat} and \eqref{eq:rhohat}, we  express $\phi_{1,\omega_0}$ as
\begin{equation}\label{def:phi_1versus-ro}
\phi_{1,\omega_0}=\frac{-\omega_0^2(\ee^{i\omega_0}-1)}{q(\omega_0)} \Delta_{0} \rho_{2,\omega_0} +\frac{i\omega_0^3(\ee^{i\omega_0}+1)}{q(\omega_0)} \Delta_{(0,0)} \rho_{1,\omega_0}=
\frac{\omega_0^2}{s(\omega_0)} \left( \sin(\omega_0/2) \Delta_{0} \rho_{2,\omega_0} - \omega_0  \cos(\omega_0/2) \Delta_{(0,0)} \rho_{1,\omega_0} \right),
\end{equation}
with $q(\omega_0)$ and $s(\omega_0)$ in (\ref{def:p&q}) and (\ref{eq:sw0}), respectively.
While either of the summands in \eqref{def:phi_1versus-ro} is only partially localized and still includes a sinusoidal trend, it is the combination of both that results in the cancelation of all residual components.
In a similar way
\begin{equation}\label{def:phi_2versus-ro}
\begin{array}{lll}
\phi_{2,\omega_0}&=&\frac{\omega_0^2(1-\ee^{i\omega_0})}{q(\omega_0)} \Delta_{0} \rho_{1,\omega_0}-
\frac{i\omega_0}{(1-\ee^{i\omega_0})q(\omega_0)} \Big(
2i\omega_0 \ee^{i\omega_0} \Delta_{(\omega_0,-\omega_0)}  \rho_{2,\omega_0} (\bullet+1)
+(1-\ee^{2i\omega_0}) \Delta_{(0,0)} \rho_{2,\omega_0}(\bullet+1)  \Big)\smallskip\\
&=&\frac{\omega_0}{s(\omega_0)} \left(  \omega_0 \sin(\omega_0/2) \Delta_{0} \rho_{1,\omega_0}  -\frac{\omega_0}{2 \sin(\omega_0/2)} \Delta_{(\omega_0,-\omega_0)}  \rho_{2,\omega_0} (\bullet+1) +\cos(\omega_0/2) \Delta_{(0,0)} \rho_{2,\omega_0} (\bullet+1) \right),
\end{array}
\end{equation}
where $q(\omega_0)$ and $s(\omega_0)$ are given in (\ref{def:p&q}) and (\ref{eq:sw0}), respectively.
Thus $\phi_{2,\omega_0}$ is localized in $[-1,1]$.\\

\smallskip \noindent Collecting the previous arguments, we now prove the following result.

\begin{proposition}The exponential spline space $S_{{\cal E}_4}^1(\Z)$ can be written as $S_{{\cal E}_4}^1(\Z) = S_{{\cal E}_4}(\Z) + S_{{\cal E}_3}(\Z)$.
\end{proposition}
\begin{pf}
We simply observe that a cardinal exponential spline for the space $S_{{\cal E}_4}(\Z)$ (see, e.g,  \cite{Micchelli1976} or \cite{Unser2005}) admits a unique expansion of the type
$$
s(x)=\sum_{n \in \Z} a[n] \rho_{1,\omega_0} (x-n) \quad \Leftrightarrow \quad \Lop_{1,\omega_0} s(x)=\sum_{n \in \Z} a[n] \delta(x-n), $$
where $a[n]$ is a sequence of slow growth.  The same holds for the space $S_{{\cal E}_3}(\Z)$ and the Green's function $\rho_{2,\omega_0}$ associated with the differential operator $\Lop_{2,\omega_0}$. This, in view of \eqref{def:ro_1versus-phi} and \eqref{def:ro_2versus-phi}, implies that $S_{{\cal E}_4}(\Z) + S_{{\cal E}_3}(\Z)\subset S_{{\cal E}_4}^1(\Z)$.
On the other hand from (\ref{def:phi_1versus-ro}) and (\ref{def:phi_2versus-ro}) we see that any function in $S_{{\cal E}_4}^1(\Z)$ is also in $S_{{\cal E}_4}(\Z) + S_{{\cal E}_3}(\Z)$, so completing the proof. \eop
\end{pf} \

\section{Riesz basis property}\label{sec:5}
\label{Sec:Riesz}
In this section we show that the system of integer translates of the Hermite exponential spline expansion in (\ref{eq:S41}) is stable. Indeed  we prove that, for the  vector function $\boldsymbol \phi_{\omega_0}=(\phi_{1,\omega_0},\phi_{2,\omega_0})^T$, there exist two constants $0<\alpha \le \beta <+\infty$ such that
$$
\alpha \|\M a\|_{\ell_2} \le \| \sum_{n \in \Z} \M a^T[n] \boldsymbol \phi_{\omega_0}(\cdot-n)\|_{L_2} \le \beta \|\M a\|_{\ell_2},\quad \hbox{with}\quad \M a\in \ell_2^{2\times 1}(\Z).
$$
The result is stated in the following theorem.

\begin{theorem} The system of (vector) functions $\{\boldsymbol \phi_{\omega_0}(\cdot-n),\ n\in \Z\}
$, $\boldsymbol \phi_{\omega_0}=(\phi_{1,\omega_0}, \phi_{2,\omega_0} )^T$ with $\phi_{j,\omega_0},\ j=1,2$
as in (\ref{eq:phi12_glob}), forms a Riesz basis. \end{theorem}

\begin{pf}
We start by computing the Hermitian Fourier Gram matrix of the basis, which is given by
\begin{align*}
\widehat{\M G}(\ee^{i \omega },\omega_0)&= \sum_{k \in \Z} \widehat {\boldsymbol \phi}_{\omega_0}(\omega+2 \pi k) \, \widehat {\boldsymbol \phi}_{\omega_0}(\omega+2 \pi k)^H \\
&=\left[
\begin{array}{ll}
\displaystyle \sum_{n \in \Z} \langle \phi_{1,\omega_0},\phi_{1,\omega_0}(\cdot-n)\rangle\, \ee^{- i \omega n}  & \displaystyle \sum_{n \in \Z} \langle \phi_{1,\omega_0},\phi_{2,\omega_0}(\cdot-n)\rangle\, \ee^{- i \omega n}
\\[2ex]
\displaystyle \sum_{n \in \Z} \langle \phi_{2,\omega_0},\phi_{1,\omega_0}(\cdot-n)\rangle\, \ee^{- i \omega n} &   \displaystyle \sum_{n \in \Z} \langle \phi_{2,\omega_0},\phi_{2,\omega_0}(\cdot-n)\rangle\, \ee^{- i \omega n} \\
\end{array}
\right] \\[1ex]
&=\left[
\begin{array}{cc}
a(\omega_0) (\ee^{-i \omega}+\ee^{i \omega}) + b(\omega_0)  &  c(\omega_0) (\ee^{-i \omega}-\ee^{i \omega}), \\
\overline{c(\omega_0)} (\ee^{i \omega}-\ee^{-i \omega}) & d(\omega_0) (\ee^{-i \omega}+\ee^{i \omega}) + e(\omega_0)
\end{array}
\right]\\[1ex]
&=\left[
\begin{array}{cc}
2 a(\omega_0) \cos(\omega) + b(\omega_0)  &  -2 c(\omega_0) i \sin(\omega), \\
2 \overline{c(\omega_0)} i \sin(\omega) & 2d(\omega_0) \cos(\omega) + e(\omega_0)
\end{array}
\right],
\end{align*}
where
\begin{eqnarray*}
a(\omega_0) &:=& \frac{\omega_0(\omega_0^2-18) \cos(\omega_0) -6 (\omega_0^2-5) \sin(\omega_0) + \omega_0(\omega_0^2-12)}{12 \omega_0 (s(\omega_0))^2}, \smallskip\\
b(\omega_0) &:=& \frac{\omega_0(\omega_0^2+3) \cos(\omega_0) -3 (\omega_0^2+5) \sin(\omega_0) + \omega_0(\omega_0^2+12) }{3 \omega_0 (s(\omega_0))^2}, \smallskip\\
c(\omega_0) &:=& \frac{5 \omega_0 (\omega_0^2+3) \cos(\omega_0/2)+ \omega_0 (\omega_0^2-15) \cos(3 \omega_0/2)-72 \sin(\omega_0/2)-6 (\omega_0^2-4) \sin(3 \omega_0/2)}{24 \omega_0^2 \sin(\omega_0/2) (s(\omega_0))^2}, \smallskip\\
d(\omega_0) &:=& \frac{6 (7 \omega_0^2+6) \sin(\omega_0)+6 (\omega_0^2-3) \sin(2 \omega_0)
-\omega_0 \left( 2(7 \omega_0^2-30) \cos(\omega_0)+(\omega_0^2-12) \cos(2 \omega_0)+3(\omega_0^2+24) \right)
}{48 \omega_0^3 \sin^2(\omega_0/2) (s(\omega_0))^2},\\
e(\omega_0) &:=& \frac{-12 (2 \omega_0^2+3) \sin(\omega_0)-3 (5 \omega_0^2-6) \sin(2 \omega_0)
+2 \omega_0 \left(2 (\omega_0^2+9) \cos(\omega_0) + (\omega_0^2-18) \cos(2 \omega_0) +6 \omega_0^2 \right)}{24 \omega_0^3 \sin^2(\omega_0/2) (s(\omega_0))^2},
\end{eqnarray*}
are real functions, $\omega_0 \in [0,\pi]$ and $s(\omega_0)$ is defined as in (\ref{eq:sw0}).
We continue by observing that the Gram matrix $\widehat{\M G}(\ee^{i \omega },\omega_0)$ is symmetric and $2\pi$-periodic and that  the Riesz basis requirement is equivalent to (see \cite{JP2001})
\begin{equation}\label{eq:Rieszeigen}
+ \infty > \beta^2 =\max_{\omega \in [0,\pi]} \lambda_{\max} (\ee^{i \omega},\omega_0) \ge \min_{\omega \in [0,\pi]} \lambda_{\min} (\ee^{i \omega},\omega_0) = \alpha^2 >0,
\end{equation}
where $\lambda_{\max}(\ee^{i \omega},\omega_0)$ and $\lambda_{\min}(\ee^{i \omega},\omega_0)$ denote the maximum and minimum eigenvalues of $\widehat{\M G}(\ee^{i \omega},\omega_0)$ at frequency $\omega$, respectively. To prove (\ref{eq:Rieszeigen}), we start by computing the trace of $\widehat{\M G}(\ee^{i \omega },\omega_0)$ (which is a real-valued function that equals the sum of the two eigenvalues) as
$${\rm tr}\big(\widehat{\M G}(\ee^{i \omega },\omega_0)\big)=2(a(\omega_0)+d(\omega_0))\cos(\omega)+b(\omega_0)+e(\omega_0).$$
Since both $a(\omega_0)+d(\omega_0)$ and $b(\omega_0)+e(\omega_0)$ are bounded real numbers, ${\rm tr}\big(\widehat{\M G}(\ee^{i \omega },\omega_0)\big)$
is bounded from above, and hence $\beta<+\infty$.
Moreover, since
both $b(\omega_0)-2a(\omega_0)$ and $e(\omega_0)-2d(\omega_0)$ are real positive numbers, we can write
$${\rm tr}\big(\widehat{\M G}(\ee^{i \omega },(\omega_0))\big)=2(a(\omega_0)+d(\omega_0))\cos(\omega)+b(\omega_0)+e(\omega_0) > (b(\omega_0)-2a(\omega_0))+(e(\omega_0)-2d(\omega_0))>0;$$
i.e., the trace is also positive, which means that $\beta$ is bounded from below.
In order to prove the existence of $\alpha>0$ such that \eqref{eq:Rieszeigen} is true, it suffices to compute ${\rm det}\big(\widehat{\M G}(\ee^{i \omega },\omega_0)\big)$ (which is the product of the eigenvalues)
and verify that it is positive and bounded away from 0. The computation of the determinant yields
$$
\begin{array}{lll}
{\rm det}\big(\widehat{\M G}(\ee^{i \omega },\omega_0)\big)&=&\Big(2a(\omega_0) \cos(\omega)+b(\omega_0)\Big) \Big(2d(\omega_0) \cos(\omega) +e(\omega_0)\Big) -4\sin^2(\omega)(c(\omega_0))^2 \\
\\
&=& A(\omega_0) \cos(2\omega)+ B(\omega_0) \cos(\omega)+C(\omega_0)\,,
\end{array}
$$
with
$$
A(\omega_0)=2\Big(a(\omega_0)d(\omega_0)+ (c(\omega_0))^2\Big), \; B(\omega_0)=2\Big(a(\omega_0)e(\omega_0)+b(\omega_0)d(\omega_0)\Big), \;
C(\omega_0)=2\Big(a(\omega_0)d(\omega_0)- (c(\omega_0))^2\Big)+b(\omega_0)e(\omega_0).$$
Next we construct the lower bound
 $${\rm det}\big(\widehat{\M G}(\ee^{i \omega },\omega_0)\big)\ge C(\omega_0)-|B(\omega_0)|-|A(\omega_0)|=C(\omega_0)+B(\omega_0)-A(\omega_0)=:G(\omega_0).$$
The final step is to observe that the auxiliary function
\begin{align*}
G(\omega_0)&=\frac{180 \omega_0 \sin(\omega_0)-9 \omega_0^3 \sin(2 \omega_0) -4 (2 \omega_0^4 - 3 \omega_0^2 - 48) \cos(\omega_0) +(\omega_0^4 -24 \omega_0^2 -3) \cos(2 \omega_0)  + 7 \omega_0^4 - 78 \omega_0^2 - 189}{24 \omega_0^4 \sin^2(\omega_0/2) (s(\omega_0))^2}\\ &\ge G(0)>0
\end{align*}
is positive and increasing for $\omega_0\in [0,\pi]$, which proves existence of the lower Riesz bound. \eop
\end{pf}

\section{Re-scaled Hermite representation}\label{sec:6}
We now specify the Hermite functions with respect to the grid $h \Z$ where $h>0$ is the sampling step. The corresponding generators $
\boldsymbol \phi_{\omega_0}^h=(\phi_{1,\omega_0}^h,\phi_{2,\omega_0}^h)^T$ are obtained from
$\boldsymbol \phi_{\omega_0}^1:=\boldsymbol \phi_{\omega_0}$ and satisfy
\begin{align}
\label{eq:phih}
\left\{
\begin{array}{ll}
\phi_{1,\omega_0}^h(x)&= \phi_{1,h\omega_0}(x/h)
\\
\phi_{2,\omega_0}^h(x)&=h\, \phi_{2,h\omega_0}(x/h),
\end{array}\right.
\end{align}
where $\phi_{j,h\omega_0},\ j=1,2$ are the Hermite cardinal functions in \eqref{eq:phi12_glob} with $\omega_0$ replaced by $h \omega_0$. Note that the second function is re-normalized to fulfill the Hermite interpolation condition $(\phi_{2,\omega_0}^h)'(0)=1$. Likewise, the derivatives  satisfy the scaling relation
\begin{align}
\label{eq:phiDh}
\left\{
\begin{array}{ll}
(\phi_{1,\omega_0}^h)'(x)&=\frac{1}{h} (\phi_{1,h\omega_0})'(x/h)
\\
(\phi_{2,\omega_0}^h)'(x)&=(\phi_{2,h\omega_0})'(x/h).\end{array}\right.
\end{align}
We then define the Hermite spline space at resolution $h$ as
\begin{align}\label{Eq:hermitespaceh}
S_{{\cal E}_4}^1(h\Z)=\left\{s_h(x)=\sum_{n \in \Z} \M a_h^T[n] {\boldsymbol \phi}^{h}_{\omega_0}(x-nh) : \M a_h \in \ell_2^{2\times 1}(\Z) \right\}.
\end{align}
The asymptotic behaviour of the re-scaled Hermite functions $\phi_{j,\omega_0}^h,\ j=1,2$ is investigated in the next proposition.

\begin{proposition}
\label{Prop:converg}
The re-scaled Hermite functions $\phi_{j,\omega_0}^h,\ j=1,2$ satisfy
$$
\lim_{h\rightarrow 0}\phi_{1,\omega_0}^h(hx)=\left\{
\begin{array}{ll}
(-2x+1)(x+1)^2, & \mbox{for } -1\le x \le0, \\
(2x+1)(x-1)^2,&  \mbox{for } \ 0<x\le 1, \end{array}\right.
$$
$$
\lim_{h\rightarrow 0}\frac{1}{h} \phi_{2,\omega_0}^h(hx)=\left\{
\begin{array}{ll}
x(x+1)^2, & \mbox{for } -1\le x \le0, \\
x(x-1)^2,&  \mbox{for } \ 0<x\le 1.
\end{array}\right.
$$
\end{proposition}
\begin{pf}
In light of \eqref{eq:phih}, the result is obtained simply by taking the limit of (\ref{def:g_real-version}) as $\omega_0\rightarrow 0$. \eop
\end{pf}

\medskip \noindent
This result is important because it shows that the re-scaled Hermite functions converge to the cardinal Hermite cubic splines as $h \to 0$.
\begin{remark}
The implication of Proposition \ref{Prop:converg} is that the asymptotic properties of the exponential Hermite splines are the same as those of the classical cubic Hermite splines. They are therefore endowed with the same fourth-order of approximation. This happens to be the order of approximation of the cubic B-splines, which are included in the space spanned by the Hermite splines as $\omega_0\to 0$.
\end{remark}

\section{Multiresolution properties}\label{sec:7}
\label{Sec:multi}
To make the multiresolution structure of these spaces apparent, we
define the Hermite spline space at resolution $h$ given
in \eqref{Eq:hermitespaceh} in terms of the Green's functions $\boldsymbol \rho_{\omega_0}=(\rho_{1,\omega_0},\rho_{2,\omega_0})^T$.
To this end, we use the convolution relation
$$
\boldsymbol \phi_{\omega_0}^h(x)= \sum_{k \in \Z} \M R_h[k] \, \boldsymbol \rho_{\omega_0}(x-h k),
$$
which is the time-domain counterpart of \eqref{Eq:FouriertoGreen} when properly rescaled to the grid $h\Z$. This allows us to show that
$$
s_h(x)=\sum_{n \in \Z} \M b_h^T[n]\, \boldsymbol \rho_{\omega_0}(x-nh),
$$
where $\M b_h^T[n]=\sum_{k \in \Z} \M a_h^T[n-k]\, \M R_h[k]= (\M a_h^T \ast \M R_h)[n]$.
Since the basis functions in this second representation do not depend on $h$,
we can infer that $S_{{\cal E}_4}^1(h\Z)\subset S_{{\cal E}_4}^1(\frac{h}{m}\Z)$ for any integer $m>1$, simply because the basis functions of the coarser space are a (subsampled) subset of ones located on the finer grid.
On the side of the Hermite generators, the corresponding two-scale relation is
\begin{align}
\boldsymbol \phi_{\omega_0}^{h}(x)
=\sum_{n \in \Z} \M H_{h \to h/m}[n] \, \boldsymbol \phi_{\omega_0}^{\frac{h}{m}} \left(x-n\frac{h}{m} \right),
\label{eq:Greenrepro}
\end{align}
with refinement mask
$$
\M H_{h \to h/m}[n]=\left[
\begin{array}{cc}
\phi_{1,\omega_0}^{h}({n\frac{h}{m}}) & (\phi_{1,\omega_0}^{h})'(n\frac{h}{m}) \nonumber \smallskip\\
\phi_{2,\omega_0}^{h}(n\frac{h}{m}) & (\phi_{2,\omega_0}^{h})'(n\frac{h}{m}) \end{array}
\right]
=\left[
\begin{array}{cc}
\phi_{1,h\omega_0}\big(\frac{n}{m}\big) & \frac{1}{h} (\phi_{1,h\omega_0})'\big(\frac{n}{m}\big)  \smallskip\\
h\, \phi_{2,h\omega_0}\big(\frac{n}{m}\big) & (\phi_{2,h\omega_0})'\big(\frac{n}{m}\big)\end{array}
\right],
$$
which follows from the application of the Hermite interpolation formula with respect to the grid $h\Z$ as well as from \eqref{eq:phih} and \eqref{eq:phiDh}.

\smallskip \noindent
As an application of this result, we write down the $m$-ary vector subdivision scheme for computing
the function
$$
s(x)=\sum_{n \in \Z} \M a_0^T[n]\,\boldsymbol \phi_{\omega_0}(x-n),
$$
as well as its first derivative, at any arbitrary fine grid with $h_J=1/m^{J}$ starting from its values at the integers.

\smallskip \noindent For readers not familiar with subdivision, we shortly recall that a vector subdivision scheme is an efficient iterative procedure based on the repeated application of refinement rules transforming, at each iteration, a sequence of vectors into a denser sequence of vectors. Whenever convergent, they generate vector functions related to the vector data used to start the iterative procedure (see \cite{DyL2002} or \cite{Warren2002} for details on subdivision schemes). The present subdivision scheme turns out to be interpolatory: since each finer sequence contains the coarser one, the initial vector data corresponds to the samples of the limit function. We refer the reader to \cite{CoCoSa07_alesund} and \cite{CZ04} for theoretical results on interpolatory vector subdivision schemes. Moreover, our vector subdivision scheme is of Hermite-type, with the understanding that the initial data and the vectors generated at each step represent function values and consecutive derivatives up to a certain order. Details on interpolatory as well as non-interpolatory Hermite subdivision schemes can be found in \cite{CMR14,ContiCotroneiSauer14,MerrienSauer12}.

\smallskip \noindent
Concretely, the interpolatory Hermite-type subdivision algorithm associated to (\ref{eq:Greenrepro}) proceeds recursively for $j=0,\dots,J-1$ by computing for all $n \in \Z$
\begin{align}
\label{eq:subsivid}
\M a_{j+1}[n]=\sum_{\ell \in \Z} \M H_{j}[mn-\ell] \, \M a_{j}[\ell],
\end{align}
where $\M H_{j}[n]:=\M H^T_{h_{j}  \to h_{j+1}}[n]$ and $h_j=\frac{1}{m^j}$.
When $m=2$ (dyadic Hermite interpolation), each step
involves an upsampling by a factor of two followed by a matrix filtering. The corresponding dyadic filters (or dyadic subdivision masks) $\{\M H_{j}[n],\ j\ge 0\}$, which are non-zero for the entries $n=-1,0,1$ only, are described by the matrix sequences
\begin{eqnarray}\label{def:maschera_binaria}
\M H_{j}[-1]&=&\left(
\begin{array}{lr}
\frac{1}{2} & \frac{1-\ee^{i \omega^{(j+1)}_0}}{2 i \omega^{(j)}_0 (\ee^{i \omega^{(j+1)}_0}+1)} \times \frac{1}{2^j} \smallskip \\
\frac{i \omega^{(j)}_0 (\ee^{i \omega^{(j+1)}_0}-1)^2}{D(\omega^{(j)}_0)} \times 2^j & \tfrac{i \omega^{(j)}_0 \ee^{i \omega^{(j+1)}_0}-\ee^{i \omega^{(j)}_0}+1}{D(\omega^{(j)}_0)}
\end{array}\right)=
\left(
\begin{array}{lr}
\frac{1}{2} & -\frac{\tan(\omega_0^{(j)}/4)}{2 \omega_0^{(j)}} \times h_j \\
\frac{2 \omega_0^{(j)} \sin^2(\omega_0^{(j)}/4)}{s(\omega_0^{(j)})} \times \frac{1}{h_j} &
\frac{2 \sin(\omega_0^{(j)}/2)-\omega_0^{(j)}}{2 s(\omega_0^{(j)})}
\end{array}\right), \medskip \nonumber\\
\M H_{j}[0]&=&\left(
\begin{array}{cc}
1 & 0 \\
0 & 1
\end{array}\right), \medskip\\
\M H_{j}[1]&=&
\left(
\begin{array}{lr}
\frac{1}{2} & -\frac{1-\ee^{i \omega^{(j+1)}_0}}{2 i \omega^{(j)}_0 (\ee^{i \omega^{(j+1)}_0}+1)} \times \frac{1}{2^j} \smallskip \\
-\frac{i \omega^{(j)}_0 (\ee^{i \omega^{(j+1)}_0}-1)^2}{D(\omega^{(j)}_0)} \times 2^j & \tfrac{i \omega^{(j)}_0 \ee^{i \omega^{(j+1)}_0}-\ee^{i \omega^{(j)}_0}+1}{D(\omega^{(j)}_0)}
\end{array}\right)=
\left(
\begin{array}{lr}
\frac{1}{2} & \frac{\tan(\omega^{(j)}_0/4)}{2 \omega^{(j)}_0} \times h_j \\
-\frac{2 \omega^{(j)}_0 \sin^2(\omega^{(j)}_0/4)}{s(\omega^{(j)}_0)} \times \frac{1}{h_j} & \frac{2 \sin(\omega^{(j)}_0/2)-\omega^{(j)}_0 }{2 s(\omega^{(j)}_0)}
\end{array}\right),
\nonumber
\end{eqnarray}
where
$$\omega^{(j)}_0:=\omega_0/2^{j}=\omega_0 h_j, \qquad D(\omega^{(j)}_0):=i \omega^{(j)}_0(1+\ee^{i \omega^{(j)}_0})+2(1-\ee^{i \omega^{(j)}_0}),\quad j\ge 0$$
and $s(\omega^{(j)}_0)$ as in (\ref{eq:sw0}).
The output of the algorithm yields the sequence $\M a_{J}^T[n]=\left(s(n/2^J), s'(n/2^J)\right)$. Note that, as the refinement masks are resolution-dependent, the scheme can be categorized as being non-stationary.

\begin{remark}
The non-stationary $j$-level subdivision mask in (\ref{def:maschera_binaria}) is such that,  for $\M D=\left(
\begin{array}{lr}
1 & 0  \\
0  & \frac{1}{2}
\end{array}\right)$,
$$
\lim_{j \rightarrow \infty} \M D^j\M H_{j}[-1] \M D^{-j}=
\left(
\begin{array}{lr}
\frac{1}{2} & -\frac{1}{8}  \smallskip \\
\frac{3}{2}  & -\frac{1}{4}
\end{array}\right), \qquad
\lim_{j \rightarrow \infty} \M D^j \M H_{j}[1] \M D^{-j}=
\left(
\begin{array}{lr}
\frac{1}{2} & \frac{1}{8} \smallskip \\
-\frac{3}{2} & -\frac{1}{4}
\end{array}\right),
$$
i.e., it is asymptotically similar to Merrien's stationary scheme based on Hermite cubic splines \cite{Merrien}.
\end{remark}

\begin{remark}
Like the non-stationary subdivision scheme in \cite{Romani10}, equation \eqref{eq:subsivid} describes a
2-point Hermite subdivision scheme reproducing ellipses.
\end{remark}

\section{Equivalent B\'{e}zier representation}\label{sec:8}

The generalized Bernstein basis functions for the space ${\cal E}_4$ with $x \in [0,1]$ are special instances of exponential B-splines with multiple knots and they have been investigated by several authors \cite{Mazure1,Mazure05}.
For the sake of completeness, we recall their definition and main properties.

\smallskip \noindent In analogy with Bernstein polynomials of degree $3$, the four Bernstein basis functions $b_{\ell, \omega_0}(x),\ \ell=0,\cdots,3$ of ${\cal E}_4$
satisfying
\begin{enumerate}\label{Bproperties}
\item[i)] \emph{symmetry}: $b_{\ell, \omega_0}(x)=b_{3-\ell,\omega_0}(1-x)$ for all $\ell=0,\cdots,3$ and $x \in [0,1]$;
\item[ii)] \emph{endpoint conditions}, listed only for $b_{0,\omega_0}$ and $b_{1,\omega_0}$:\\
$$b_{0,\omega_0}(0)=1,\ b_{0,\omega_0}(1)=0,\ (b_{0,\omega_0})'(1)=(b_{0,\omega_0})^{''}(1)=0,\qquad  b_{1, \omega_0}(0)=b_{1,\omega_0}(1)=0,\ (b_{1, \omega_0})'(1)=0\,;$$
\item[iii)] \emph{partition of unity}: $\displaystyle{\sum_{\ell=0}^3 b_{\ell, \omega_0}(x)=1}$ for all $x\in[0,1]$;
\item[iv)] \emph{non-negativity}: $b_{\ell, \omega_0}(x)\ge 0$ for all $x\in[0,1]$  and $\ell=0,\cdots,3$;
\end{enumerate}
\medskip
\noindent are given by
\medskip
\begin{equation}
\begin{array}{lll}
b_{0,\omega_0}(x)&=& \frac{2 i\omega_0 \ee^{i\omega_0}}{r(\omega_0)}-\frac{2 i\omega_0 \ee^{i\omega_0}}{r(\omega_0)}x+
\frac{1}{r(\omega_0)}\ee^{i\omega_0\, x}-\frac{\ee^{2i\omega_0}}{r(\omega_0)}\ee^{-\omega_0\, x} \smallskip\\
&=&\frac{\omega_0}{\omega_0 - \sin(\omega_0)}(1-x) -\frac{\sin(\omega_0(1-x))}{\omega_0 - \sin(\omega_0)}, \medskip \\
b_{1,\omega_0}(x)&=&\frac{(1-\ee^{i\omega_0}) p(\omega_0)}{q(\omega_0) r(\omega_0)}+
\frac{i\omega_0 (\ee^{i\omega_0}-1)^3}{q(\omega_0) r(\omega_0)}x+\frac{r(\omega_0)-q(\omega_0)}{q(\omega_0) r(\omega_0)}\ee^{i\omega_0\, x}
+\frac{\ee^{i\omega_0} (p(\omega_0)-q(\omega_0)) }{q(\omega_0) r(\omega_0)}\ee^{-i\omega_0\, x}\smallskip\\
&=&
\frac{\sin(\omega_0/2)}{s(\omega_0)} -\frac{2 \omega_0 \sin^3(\omega_0/2)}{s(\omega_0) \, (\omega_0 - \sin(\omega_0))} (1-x) +
\left( \frac{1}{\omega_0 - \sin(\omega_0)} + \frac{\cos(\omega_0/2)}{s(\omega_0)} \right) \sin(\omega_0 (1-x)) -\frac{ \sin(\omega_0/2)}{s(\omega_0)} \cos(\omega_0 (1-x)), \medskip\\
b_{2,\omega_0}(x)&=&\frac{(1-\ee^{i\omega_0})}{q(\omega_0)}+\frac{i\omega_0 (1-\ee^{i\omega_0})^3}{q(\omega_0) r(\omega_0)}x+\frac{p(\omega_0)-q(\omega_0)}{q(\omega_0) r(\omega_0)}\ee^{i\omega_0\, x}+\frac{ \ee^{i\omega_0}(r(\omega_0)-q(\omega_0))  }{q(\omega_0) r(\omega_0)}\ee^{-i\omega_0\, x}\smallskip\\
&=& \frac{\sin(\omega_0/2)}{s(\omega_0)} -\frac{2 \omega_0 \sin^3(\omega_0/2)}{s(\omega_0) \, (\omega_0 - \sin(\omega_0))} x +
\left( \frac{1}{\omega_0 - \sin(\omega_0)} + \frac{\cos(\omega_0/2)}{s(\omega_0)} \right) \sin(\omega_0 x) -\frac{ \sin(\omega_0/2)}{s(\omega_0)} \cos(\omega_0 x), \medskip \\
b_{3,\omega_0}(x)&=&\frac{2 i\omega_0 \ee^{i\omega_0}}{r(\omega_0)}x+
-\frac{\ee^{i\omega_0}}{r(\omega_0)}\ee^{i\omega_0\, x}+
\frac{\ee^{i\omega_0}}{r(\omega_0)}\ee^{-i\omega_0\, x}\smallskip\\
&=&\frac{\omega_0}{\omega_0 - \sin(\omega_0)} x -\frac{\sin(\omega_0 x)}{\omega_0 - \sin(\omega_0)},
\end{array}
\end{equation}
\medskip
where  $q(\omega_0)$ is given in (\ref{def:p&q}) and $r(\omega_0):=1+2i \omega_0 \ee^{i\omega_0}-\ee^{2i\omega_0}$.

\medskip \noindent For later use, we mention that, by symmetry,
$
(b_{2,\omega_0})'(0)=(b_{3,\omega_0})'(0)=0
$. Similarly,
$$(b_{0,\omega_0})'(0)=(b_{2,\omega_0})'(1)=-(b_{1,\omega_0})'(0)=-(b_{3,\omega_0})'(1)=\frac{p(\omega_0)-r(\omega_0)}{r(\omega_0)}=\frac{\omega_0 (\cos(\omega_0) - 1)}{\omega_0 - \sin(\omega_0)}.
$$

\begin{figure}[h!]
\centering
\includegraphics[width=5.5cm]{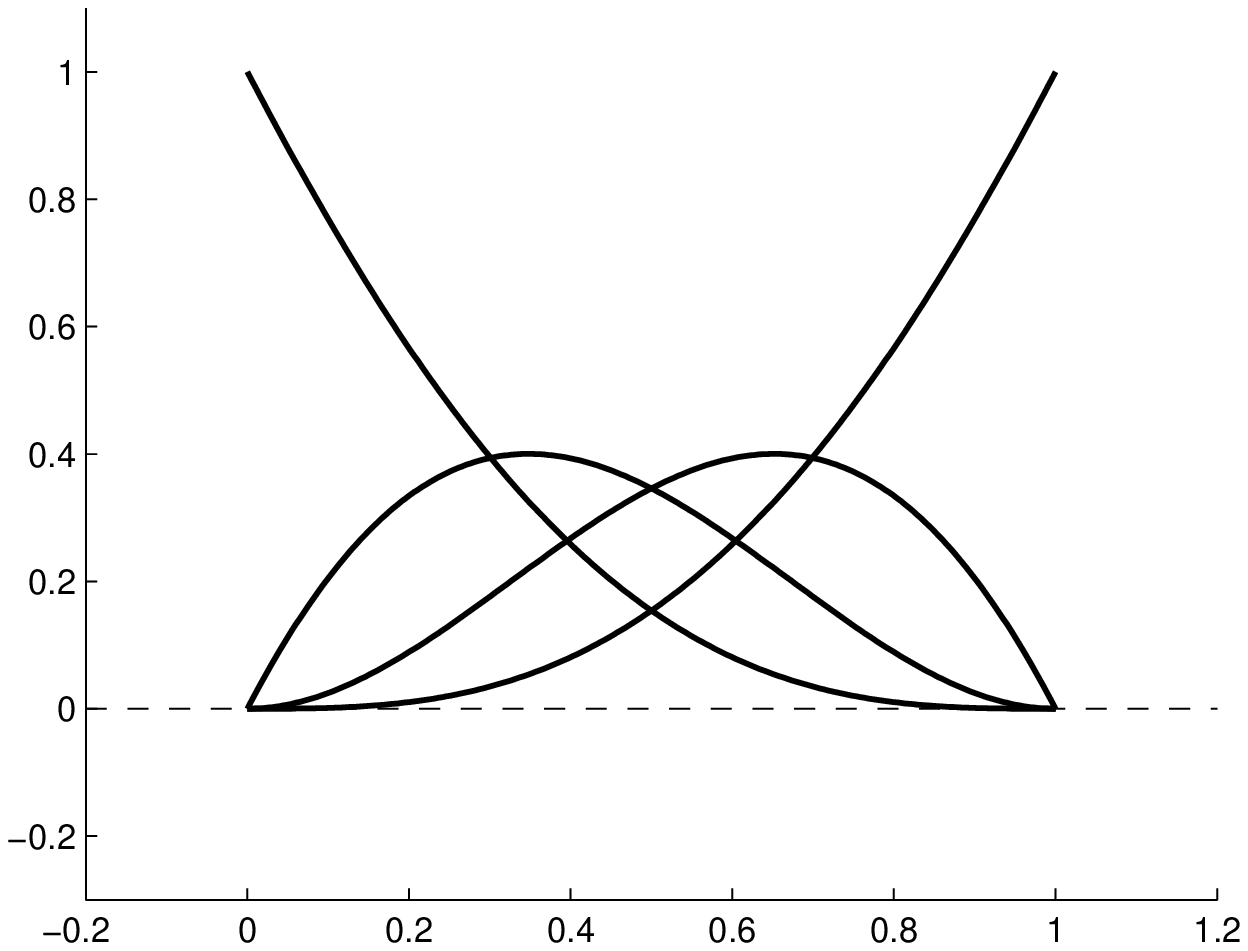}\hspace{0.5cm}
\includegraphics[width=5.5cm]{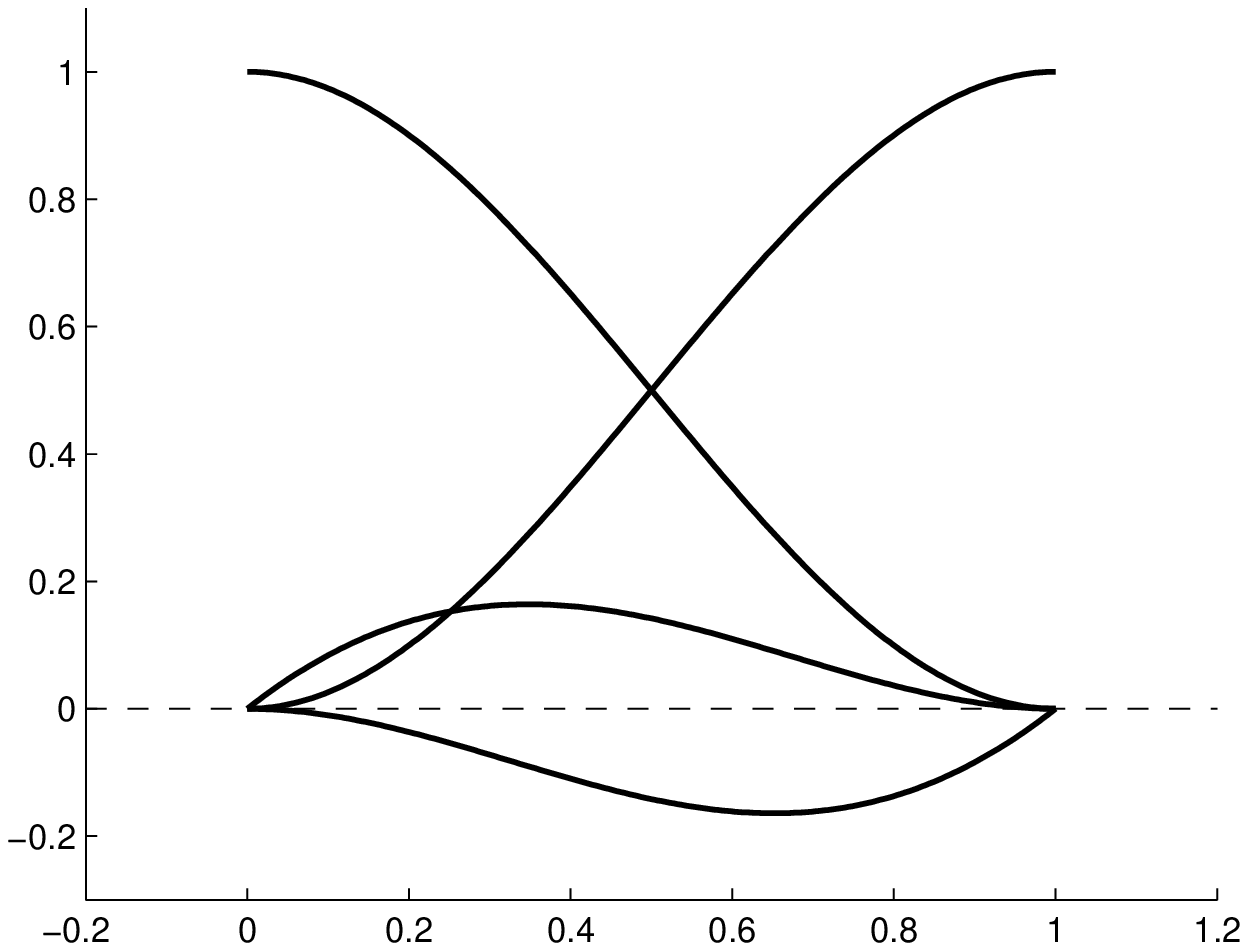}
\caption{\label{Fig:Bernstein}
Exponential Bernstein basis functions (left) versus exponential Hermite basis functions (right)
for  $\omega_0=3/4\pi$.}
\end{figure}

\noindent It is well known that cubic Hermite interpolation can be expressed in terms of cubic B\'ezier basis functions. To achieve the same in the present context,
let us consider the task of computing $b_{\ell,\omega_0}(x),\ \ell=0,\cdots, 3$ as specified by \eqref{eq:S41}
for   $x\in[n,n+1)$.
Defining $t=x-n \in[0,1)$, we simplify the expansion as
\begin{align}
\label{eq:evalhermite}
b_{\ell,\omega_0}(n+t)=b_{\ell,\omega_0}(n)\phi_{1,\omega_0}(t)+(b_{\ell,\omega_0})'(n)\phi_{2,\omega_0}(t)+b_{\ell,\omega_0}(n+1)\phi_{1,\omega_0}(t-1)+(b_{\ell,\omega_0})'(n+1)\phi_{2,\omega_0}(t-1),
\end{align}
by retaining only the four Hermite basis functions that are non-vanishing within the interval (see Figure \ref{Fig:Bernstein}).

From the endpoint condition (ii),  we readily obtain the conversion between the two types of representations as
\begin{align}
\label{eq:hermtobez}
\left[
\begin{array}{l}
\phi_{1,\omega_0}(t) \\
\phi_{2,\omega_0}(t) \\
\phi_{1,\omega_0}(t-1) \\
\phi_{2,\omega_0}(t-1)
\end{array}
\right]
=\left[
\begin{array}{cccc}
1 & 1 & 0 & 0 \\
0 & \frac{r(\omega_0)}{r(\omega_0)-p(\omega_0)} & 0 & 0 \\
0 & 0 & 1 & 1 \\
0 & 0 & -\frac{r(\omega_0)}{r(\omega_0)-p(\omega_0)} & 0
\end{array}
\right] \, \left[
\begin{array}{l}
b_{0,\omega_0}(t) \\
b_{1,\omega_0}(t) \\
b_{2,\omega_0}(t) \\
b_{3,\omega_0}(t) \\
\end{array}
\right]\,.
\end{align}

\begin{remark}
Note that $\lim_{\omega_0\rightarrow 0}\frac{r(\omega_0)}{r(\omega_0)-p(\omega_0)}=\frac13$. This indicates that the above conversion matrix provides in the limit the conversion matrix for cubic polynomial Hermite splines, as expected.
\end{remark}

\section{Link with scalar subdivision}\label{sec:9}
\medskip \noindent We conclude the paper by showing that the Hermite subdivision scheme discussed in Section \ref{Sec:multi} can also be converted into a non-uniform, non-stationary scalar subdivision scheme for exponential B-splines with double knots spanning $S^1_{{\cal E}_4}(\Z)$. This is the (new) exponential counterpart of the subdivision scheme for cubic B-splines with double knots considered in \cite{Juttler}. Based on the conversion between Hermite and B\'ezier functions for
${\cal E}_4$ given by \eqref{eq:hermtobez}, we see that, for $j\ge 0$, in the interval $[\frac{\ell}{2^j}, \frac{\ell+1}{2^j}]$, the function
\begin{equation}\label{eq:funzione}
f_j[n]\phi_{1,\omega_0}(x)+d_j[n]\phi_{2,\omega_0}(x)+f_j[n+1]\phi_{1,\omega_0}(x-1/2^j)+d_j[n+1]\phi_{2,\omega_0}(x-1/2^j)
\end{equation}
can be written as
$$
\begin{array}{ll}
f_j[n]b_{0,\omega_0}(x)+\left(f_j[n]+\frac{r(\omega_0)}{2^j(r(\omega_0)-p(\omega_0))}d_j[n]\right)b_{1,\omega_0}(x)+\\
\left(f_j[n+1]-\frac{r(\omega_0)}{2^j(r(\omega_0)-p(\omega_0))}d_j[n+1]\right)b_{2,\omega_0}(x)+f_j[n+1]b_{3,\omega_0}(x)\,,
\end{array}
$$
or, in a more compact form, as
$$
f_j[n]b_{0,\omega_0}(x)+p_j[2n+1]b_{1,\omega_0}(x)+
p_j[2n+2]b_{2,\omega_0}(x)+f_j[n+1]b_{3,\omega_0}(x)\,,
$$
with
\begin{equation}\label{eq:relazione}
\underbrace{
\left(
\begin{array}{l}
p_j[2n] \\
p_j[2n+1]
\end{array}\right)}_{\M p_{j}[n]}=\underbrace{\left(
\begin{array}{cc}
1 & -\frac{r(\omega_0)}{2^j(r(\omega_0)-p(\omega_0))}\\
1 & \frac{r(\omega_0)}{2^j(r(\omega_0)-p(\omega_0))}
\end{array}\right)}_{\M M_{j}} \,  \underbrace{ \left(
\begin{array}{l}
f_j[n] \\
d_j[n]
\end{array}\right)}_{\M a_j[n]},\quad n\ge 0\,.
\end{equation}

\begin{figure}[h!]
\centering
\includegraphics[width=10cm]{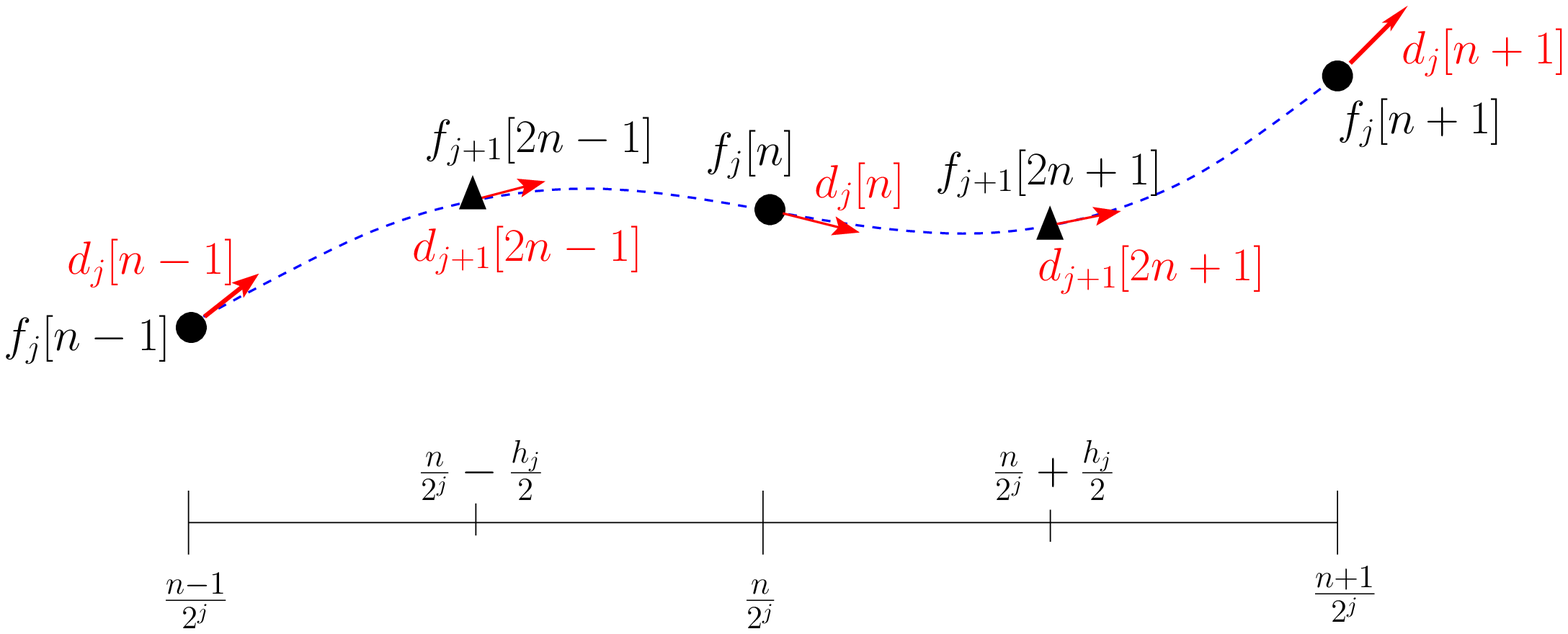} \\\vspace{0.4cm}
\hspace{-0.4cm}\includegraphics[width=7.75cm]{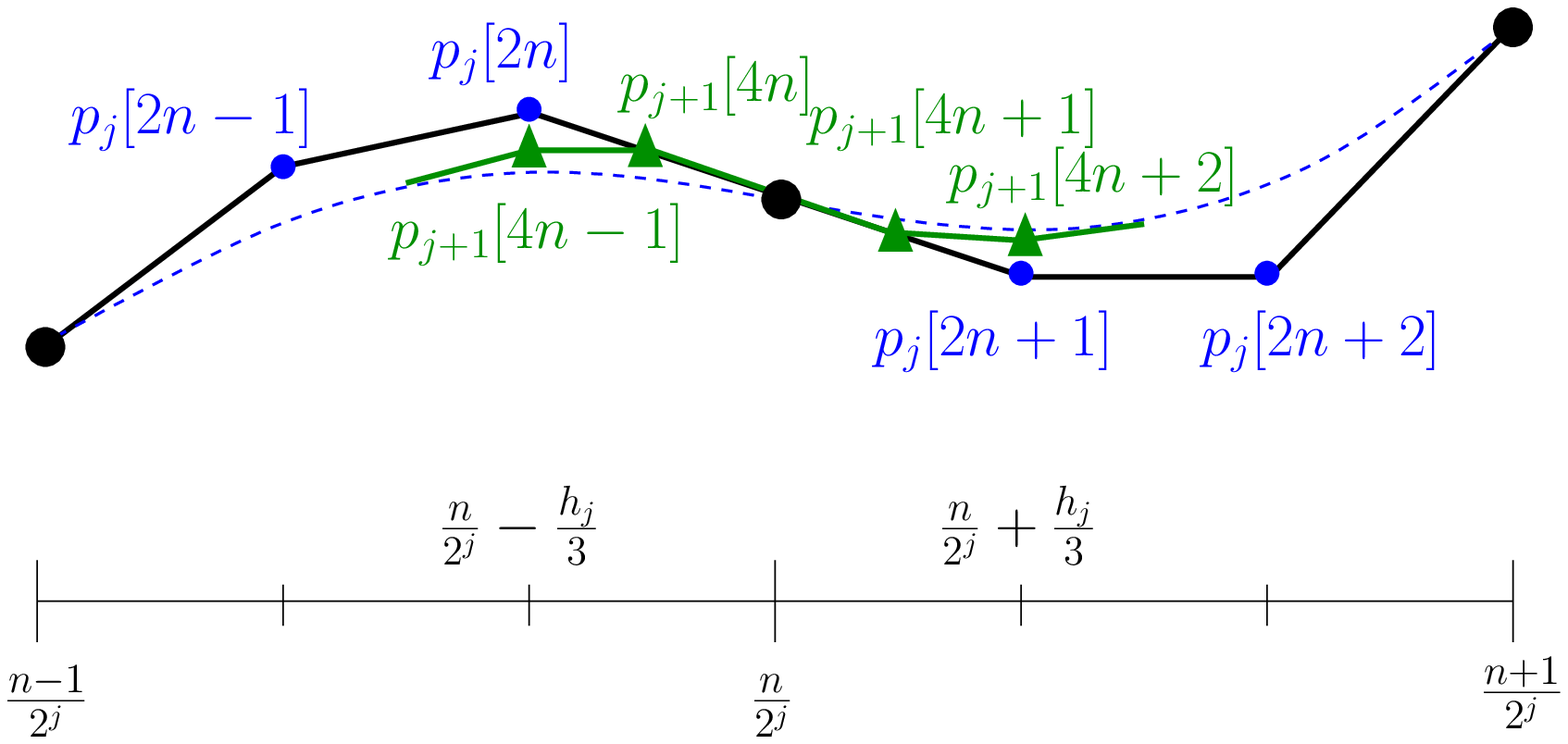}
\caption{Geometric interpretation of the subdivision schemes in \eqref{eq:2-point_Hermite} and \eqref{def:4-point}.
Top: at each step the interpolatory vector subdivision scheme \eqref{eq:2-point_Hermite} creates a new vector between any two old vectors and retains them.
Bottom: at each step the approximating scalar subdivision scheme \eqref{def:4-point} creates two new control points between any two old ones and discards them.}
\label{fig:four-point}
\end{figure}

\noindent
At this point, we recall that the dyadic Hermite subdivision scheme with mask \eqref{def:maschera_binaria} and the repeated evaluation of the local Hermite interpolant at interval mid points (see, for example, \cite{Merrien}) can be explicitly written as
\begin{equation}\label{eq:2-point_Hermite}
\M a_{j+1}[2n]:=\M a_{j}[n],\qquad \M a_{j+1}[2n+1]:=\M H_j[1]\, \M a_{j}[n]+\M H_j[-1]\M a_{j}[n+1]\,,
\end{equation}
or, in view of \eqref{eq:relazione}, as
\begin{equation}\label{def:4-point}
\M p_{j+1}[2n]:=\M M_{j+1}\, \M M^{-1}_{j}\, \M p_{j}[n],\qquad \M p_{j+1}[2n+1]:=\M M_{j+1}\,\M H_j[1]\, \M M^{-1}_{j}\, \M p_{j}[n]+\M M_{j+1}\,\M H_j[-1]\, \M M^{-1}_{j}\, \M p_{j}[n+1]\,.
\end{equation}
Since at each iteration $j$, the latter formulas define $p_{j+1}[4n],\ p_{j+1}[4n+1] ,\ p_{j+1}[4n+2],\ p_{j+1}[4n+3]$, the vector rules in (\ref{def:4-point}) do identify four scalar rules that we can associate to a non-uniform and non-stationary scalar  subdivision scheme. This is the exponential counterpart of the four-point scheme in \cite{Juttler}, whose geometric meaning is shown in Figure \ref{fig:four-point}.

\section*{Acknowledgements}
We thank Virginie Uhlmann for her helpful comments on the manuscript, and for her development of a companion algorithm for image segmentation whose preliminary results were presented in \cite{Noi2014}.\\
Support from the Italian GNCS-INdAM and the Swiss Science Foundation under Grant 200020-144355 is gratefully acknowledged. Lucia Romani also acknowledges the support of MIUR-PRIN 2012 (grant 2012MTE38N).


\end{document}